\documentclass[12pt,english]{amsart}

\usepackage{amscd}
\usepackage[sc,osf,slantedGreek]{mathpazo}
\usepackage{euscript}
\usepackage{graphics}
\usepackage{color}
\usepackage[all]{xy}
\usepackage[T1]{fontenc}

\newcommand{\B}[1]{{\mathbb #1}}
\newcommand{\C}[1]{{\EuScript #1}}

\newtheorem{theorem}[subsection]{Theorem}%[section]
\newtheorem{corollary}[subsection]{Corollary}
\newtheorem{lemma}[subsection]{Lemma}
\newtheorem{proposition}[subsection]{Proposition}

\theoremstyle{definition}

\newtheorem{example}[subsection]{Example}
\theoremstyle{remark}
\newtheorem{remark}[subsection]{Remark}

\numberwithin{figure}{section}
\numberwithin{table}{section}
\newcommand{\al}{{\alpha}}

\newcommand{\be}{{\beta}}

\newcommand{\Om}{{\Omega}}
\newcommand{\om}{{\omega}}

\newcommand{\ga}{{\gamma}}

\newcommand{\la}{{\lambda}}

\newcommand{\Mo}{(M,\omega )}

\newcommand\Diff{\operatorname{Diff}}
\newcommand\Flux{\operatorname{Flux}}
\newcommand\Symp{\operatorname{Symp}}
\newcommand\Ham{\operatorname{Ham}}

\newcommand\Hom{\operatorname{Hom}}

\begin{document}

\title{A cocycle on the group of symplectic diffeomorphisms}
\author{\'Swiatos\l aw R. Gal}
\address{Instytut Matematyczny,
Uniwersytet Wroc\l awski,
pl.~Grunwaldzki 2/4,
Wroc\l aw,
Poland}
\address{Department of Mathematics,
The Ohio State University,
100 Math Tower,
231 West 18th Avenue,
Columbus, {\sc oh} 43210-1174, {\sc usa}}
\email{sgal@math.uni.wroc.pl}
\urladdr{http://www.math.uni.wroc.pl/\~{}sgal}
\author{Jarek K\k edra}
\curraddr{Mathematical Sciences,
University of Aberdeen, 
Meston Building, 
Aberdeen {\sc ab243ue}, 
Scotland, {\sc uk}}
\curraddr{Institute of Mathematics, University of Szczecin,
ul.~Wielkopolska 15, 70-451 Szczecin, Poland}
\email{kedra@maths.abdn.ac.uk}
\urladdr{http://www.maths.abdn.ac.uk/\~{}kedra}
\date{\today}
\subjclass{Primary {\sc 57s25}; Secondary {\sc 20f69}}
\begin{abstract}
We define a cocycle on the group of Hamiltonian diffeomorphisms
of a symplectically aspherical manifold and investigate its properties.
The main application is an alternative proof
of the Polterovich theorem about the distortion of  cyclic
subgroups in finitely generated groups of Hamiltonian
diffeomorphisms.

\setcounter{tocdepth}{1}
\tableofcontents
\end{abstract}

\maketitle

\section{Introduction}\label{S:intro}

%Let $\Mo$ be a symplectic manifold and let $\Gamma $ be a discrete group.
In the present paper we investigate the properties of a certain
cocycle defined on the group of Hamiltonian diffeomorphisms of 
a symplectically aspherical manifold. 
The cocycle appears more or less explicitly in Gambaudo-Ghys \cite{MR1452855}
and Arnold-Khesin \cite[p.~247]{MR1612569} for the case of a symplectic ball.

As an application we proof a theorem due to Polterovich
\cite{MR2003i:53126}
that a finitely generated subgroup of the group of
Hamiltonian diffeomorphisms of a symplectically hyperbolic
manifold has undistorted cyclic subgroups.

\subsection{Distortion of cyclic subgroups}\label{SS:undis}
Let $\Gamma $ be a finitely generated group and let $|g|$ denote
the word norm of $g\in \Gamma $ with respect to some fixed finite
set of generators. A cyclic subgroup $\left<g\right>\subset \Gamma $
generated by $g$ is called {\bf undistorted} if 
there exists a positive constant $C$ such that for any integer 
$n\in \B Z$ there holds $|g^n|\geq C\cdot |n|.$
\begin{remark}
According to Gromov's definition \cite{MR1253544} finite
subgroup are also undistorted. It follows from our
definition that an undistorted cyclic subgroup is infinite.
\end{remark}

\begin{example}[Groups with undistorted cyclic subgroups]\label{E:undist}
\begin{enumerate}
\item
Free and free Abelian groups.
\item
Hyperbolic groups. 
\item
Groups acting properly discontinuously on contractible $\operatorname{CAT}(0)$-spaces.
\end{enumerate}
\end{example}

\begin{example}[Groups with distorted cyclic subgroups]\label{E:dist}
\begin{enumerate}
\item
Groups with torsion.
\item
Heisenberg groups of the uppertriangular integer matrices.
\item
Non-uniform irreducible lattices in semisimple Lie groups 
(see Lubotzky-Mozes-Raghunathan \cite{MR1828742}).
\end{enumerate}
\end{example}

The first application is the following result which was
also proved by Frauenfelder and Schlenk
in \cite[Theorem 3.7]{FS} by a different method.

\begin{theorem}\label{T:undistorted}
Let $\Mo$ be an exact symplectic manifold convex at infinity.
Let $h\colon \Gamma \hookrightarrow \Ham_c(M,\om)$ be an effective Hamiltonian
action of a finitely generated group. Then any non-trivial cyclic
subgroup of $\Gamma $ is undistorted with respect to 
the word metric given by a finite set of generators.
\end{theorem}

\begin{remark}
It follows from the Thurston Stability (see for example Theorem 3.4 in
Franks \cite{MR2288284}) that a group acting smoothly with compact
supports on an open manifold has many undistorted cyclic subgroups.
We show that the reasons for the undistortion are not purely
topological.  Namely, in Section \ref{SS:heisenberg}, we present an
example, due to \L ukasz Grabowski, of a smooth compactly supported
action of the Heisenberg group on $\B R^n$ for $n\geq 2$.  It is,
however, not clear if the above theorem holds for volume preserving
actions.
\end{remark}

A symplectic manifold is called {\bf exact} if the symplectic form
is exact. Notice that in this case $M$ is not closed. The boundary $\partial M$ 
of a symplectic manifold $\Mo$ is called {\bf convex} if there exists a vector
field $X$ defined in a neighbourhood of $\partial M$ such that 
$\C L_X\om=\om$ and it points outwards. A symplectic manifold $\Mo$ is
{\bf convex at infinity} if there exists an increasing sequence of
compact submanifolds with convex boundaries exhausting $M$.

\begin{example}[Symplectic manifolds convex at infinity]\label{E:convex}
\begin{enumerate}
\item
The standard symplectic Euclidean space $(\B R^{2n},\om)$.
\item
The cotangent bundle $(T^*X,d\la)$ over a smooth manifold $X.$
\item
Stein manifolds.
\end{enumerate}
\end{example}

Roughly speaking our proof of Theorem \ref{T:undistorted} consists of
the following steps.  We define a cocycle $\C K\colon \Ham_c\Mo \to
C_c^{\infty}(M)$ and compose it with a norm given by the
oscillation. This induces a norm on the group $\Ham_c\Mo$ and we prove
that every nontrivial cyclic subgroup of $\Ham_c\Mo$ is undistorted
with respect to this norm. Finally we prove that the action $h\colon
\Gamma \hookrightarrow \Ham_c(M,\om)$ is Lipschitz with respect to the
word norm on $\Gamma $ and the above mentioned norm on $\Ham_c\Mo.$
The proof is presented in Section~\ref{SS:undistorted}.

Let $p\colon\widetilde M\to M$ be the universal covering.
Making an observation (Proposition \ref{P:eval})
that there exists a lifting homomorphism
$\Psi\colon\Ham\Mo\to \Ham(\widetilde M,p^*\om),$ we extend the
above argument to prove (in Section \ref{S:osc}) the following 
result by Polterovich.

\begin{theorem}[Polterovich, 1.6.C in \cite{MR2003i:53126}]\label{T:polterovich}
Let $\Mo$ be a closed symplectically hyperbolic manifold.
Let $\Gamma \hookrightarrow \Ham\Mo$ be an effective Hamiltonian action
of a finitely generated group. Then $\Gamma $ has undistorted cyclic subgroups.
\end{theorem}
\noindent
Both proofs rely on the existence of certain fixed points proved 
by Schwarz \cite{MR1755825} for closed manifolds and Frauenfelder-Schlenk
\cite {MR2342472} for manifolds convex at infinity.

%\begin{definition}

A symplectic manifold $(M,\omega)$ is said to be
{\bf symplectically aspherical} if the pull-back
$p^*\omega$ of the symplectic form $\omega$
to the universal cover 
is exact, ie. $p^*\omega=d\alpha$ for some one-form ~$\alpha$.
%\end{definition}
Observe that the symplectic asphericity is equivalent to the
vanishing of the symplectic form on spheres. That is
$\int_{S^2}f^*\om=0$ for any smooth map
$f\colon S^2\to M.$ Clearly, the symplectic asphericity
depends only on the cohomology class of the symplectic form.

\begin{example}[Closed symplectically aspherical manifolds]\label{E:sas}
Aspherical symplectic manifolds, such as
tori, surfaces of positive genus and their products, are obviously 
symplectically aspherical. 
A bran\-ched covering
of a symplectically aspherical manifold is again symplectically
aspherical. This allows to produce more interesting examples
possibly with non-trivial second and higher homotopy groups
as it was done by Gompf in \cite{MR99k:57068}. Also, a symplectic
submanifold of a symplectically aspherical manifold is symplectically
aspherical (see \cite{MR2103543,MR2308893} for other constructions 
and properties).
\end{example}

%\begin{definition}\label{D:hyperbolic}
Let $\Mo$ be a closed symplectic manifold equipped with a Riemannian
metric $g$ and let $p\colon \widetilde M\to M$ be its universal cover.
The manifold $\Mo$ is called {\bf symplectically hyperbolic} if it is
symplectically aspherical and there exists a primitive $\alpha$ of
$p^*\om$ which is bounded with respect to the metric induced by $g$.
That is there exists a constant $C\in \B R$ such that
$$
\sup_{x\in \widetilde M}\|\alpha(x)\|<C,
$$
where
$\|\al(x)\|:= \sup\{\al(x)(X)\,|\,X\in T_xM \text{ and } |X|\leq 1\}.$
%\end{definition}

Notice that the symplectic hyperbolicity does not depend on the choice
of a Riemannian metric. Also, it depends only on the cohomology class of
the symplectic form. Closed oriented surfaces of genus at least 2 
and their products are 
symplectically hyperbolic as well as their branched covers and submanifolds.
A symplectically aspherical manifold with word hyperbolic fundamental group
is symplectically hyperbolic \cite[Corollary 1.13]{hyperbolic}.

\subsection*{Acknowledgements}
We thank Michael Weiss and the University of Aberdeen for the support of the
visit of the first author in Aberdeen in May 2008.
The second author thanks Fr\'ed\'eric Bourgeois, Felix Schlenk and Thomas Vogel 
for discussions during the symplectic geometry seminar in Brussels in February 2008.

%%%%%%%%%%%%%%%%%%%%%%%%%%%%%%%%%%%%%%%%%%%%%%%%%%%%%%%%%%%%%%%
\section{A cocycle on $\Symp(M,\om=d\alpha)$.}\label{S:K}
%%%%%%%%%%%%%%%%%%%%%%%%%%%%%%%%%%%%%%%%%%%%%%%%%%%%%%%%%%%%%%%

\subsection{Definition of the cocycle}\label{SS:def}
Let $C^{\infty}(M)$ denote the vector space of real valued functions on $M$. 
The group $\Symp(M,\om)$ acts on functions from the right by the 
composition. We assume that $M$ is connected throughout the paper.

Let $\alpha \in \Om^1(M)$ be a primitive of the symplectic
form $\om$. That is $d\al=\om$. Let $f\in \Symp(M,\om)$
be a symplectic diffeomorphism. 
The one-form $f^*\al - \al$ is closed. 
%Indeed,

%$$d(f^*\al - \al) = f^*d\al - d\al = f^*\om - \om = 0.$$

\begin{proposition}\label{P:H1}
Let $f\in\Symp_0(M,d\alpha)$. The form $f^*(\alpha)-\alpha$
is exact if and only if $f$ is Hamiltonian.
\end{proposition}

\begin{proof}
A diffeomorphism $f$ is Hamiltonian if and only if
there exists a symplectic isotopy $f_t$ from the identity to $f$ such that
$$
\int_{f_t(\ga)}d\al = 0
$$
for any loop $\ga\colon S^1\to M$
(McDuff-Salamon \cite[Theorem 10.12]{MR2000g:53098}).
Then we have that
$$
0=\int_{f_t(\ga)}d\al=
 \int_{f(\ga)-\ga}\al =
\int_{\ga}f^*\al - \al 
$$
which is equivalent to the exactness of $f^*\al - \al.$
\end{proof}

Let $\widehat \Ham(M,d\alpha)\subset \Symp(M,d\alpha)$ be the
kernel of the cocycle
$\Symp(M,d\alpha)\to H^1(M;\B R)$  defined by
$f\mapsto [f^*\alpha-\alpha].$
Its connected component of the identity is the group of 
Hamiltonian diffeomorphisms. Clearly, if $H^1(M,R)=0$ then 
$\widehat \Ham(M,d\al)=\Symp(M,d\alpha)$. Enlargements of the group
of Hamiltonian diffeomorphisms of a closed symplectic manifolds
are investigated by McDuff in \cite{MR2235852}.

If $f\in \widehat \Ham(M,d\al)$ then there exists a function $\C K_{\al}(f)$ 
such that $d\C K_{\al}(f)=f^*\al - \al$. Such a function
is unique up to a constant.

%We normalize $\C K(f)$ so that it is compactly supported.

\begin{proposition}\label{P:K}
Let $(M,d\alpha)$ be an exact symplectic manifold.
The map $\C K_{\al}\colon \widehat\Ham(M,\om)\to C^{\infty}(M)/\B R$ 
defined above is a cocycle.
\end{proposition}

\begin{proof}
We have to check that $\C K_{\al}(f\circ g) =\C K_{\al}(f)\circ g +\C K_{\al}(g)$
up to a constant.
\begin{eqnarray*}
d\C K_{\al}(f\circ g)&=& (f\circ g)^*\al -\al \\
&=& g^*(f^*\al) - g^*\al + g^*\al - \al \\
&=& g^*(d\C K_{\al}(f)) + d\C K_{\al}(g)\\
&=& d(\C K_{\al}(f)\circ g)+ d(\C K_{\al}(g))\\
\end{eqnarray*}

We get that $d(\C K_{\al}(f\circ g) - \C K_{\al}(f)\circ g - \C K_{\al}(g))= 0$
which means that the function 
$\C K_{\al}(f\circ g) - \C K_{\al}(f)\circ g - \C K_{\al}(g)$
is constant.
\end{proof}

\subsection{Basic properties of $\C K_{\al}$}\label{SS:properties}
\begin{proposition}\label{P:cohomology}
If $H^1(M;\B R)=0$ then
the cohomology class $[\C K_{\al}]\in H^1(\Symp(M,d\al);C^{\infty}(M)/\B R)$
does not depend on the choice of a primitive.
\end{proposition}

\begin{proof}
Let $\beta \in \Om^1(M)$ be another primitive of $\om$. Then
$\al - \beta = dF$. A straightforward calculation shows that
$\C K_{\al}(f) - \C K_{\beta}(f) = F\circ f - F,$ which means that
$\C K_{\al}$ and $\C K_{\beta}$ are cohomologous. 
\end{proof}

\begin{proposition}\label{P:H1nonzero}
There exists an  injective homomorphism
$$\iota\colon H^1(M;\B R)\to H^1(\Ham(M,d\al);C^{\infty}(M)/\B R).$$
If $\al$ and $\beta$ are two primitives of the symplectic form
then $[\C K_{\al}]-[\C K_{\be}]=\iota[\alpha - \beta].$
\end{proposition}

\begin{proof}
Let $f\in \Ham(M,d\al)$ and let $[a]\in H^1(M;\B R)$.
Since $f$ is isotopic to the identity the closed one-form
$f^*a-a$ is exact. Define a cocycle 
$\iota(a)\colon\Ham(M,d\al)\to C^{\infty}(M)/\B R$ by
$$d(\iota(a)(f)):=f^*a - a.$$
It is straightforward to check that
this defines an injective  homomorphism 
$$\iota\colon H^1(M;\B R)\to H^1(\Ham(M,d\al);C^{\infty}(M)/\B R).$$
Now we have that
$$
d(\C K_{\al}(f) - \C K_{\be}(f))=f^*(\al-\be) - (\al-\be)
$$
which means that
$[\C K_{\al}]-[\C K_{\be}]=\iota[\alpha - \beta]$
as claimed.
\end{proof}

\begin{remark}\label{R:H1}
In the definition of $\iota$ we have only used the fact that
$f$ acts trivially on the first cohomology of $M$. Thus
$\iota$ factors through
$H^1(\Diff(M,H^1);C^{\infty}(M)/\B R)$ 
the cohomology
of the group of diffeomorphisms acting trivially on $H^1(M;\B R).$
\end{remark}

\begin{proposition}\label{P:action}
Let $f\in \Ham(M,d\al)$ be a Hamiltonian diffeomorphism generated
by the Hamiltonian function $F_t\colon M\to \B R$. Let $x,y\in M$. Then
{\small
$$
\C K_{\al}(f)(x) - \C K_{\al}(f)(y) = 
\int_{f_t(x)}\al + \int_0^1F_t(f_t(x))dt
- \left(\int_{f_t(y)} \al + \int_0^1F_t(f_t(y))dt\right),
$$
}
where $f_t$ is the isotopy from the identity to $f$ generated by $F_t$.
%In other words, there exists a constant $C\in \B R$ depending on
%$(M,d\al)$ such that 
%$$
%K_{\al}(f)(x) = \int_{f_t(x)}\al + \int_0^1F_t(f_t(x))dt + C.
%$$
\end{proposition}

\begin{proof}
Let $\ga\colon [0,1]\to M$ be a curve from $x$ to $y$ and let $D$ be a
singular disc with boundary $\partial D= \ga + f_t(y) - f(\ga) - f_t(x)$.
\begin{eqnarray*}
\C K_{\al}(f)(x) &-& \C K_{\al}(f)(y)\\ 
&=& - \int_{\ga} d\C K_{\al}\\
&=& - \int _{\ga} f^*\al - \al = \int_{\ga} \al - \int_{f(\ga)}\al \\
&=& \int_D d\al - \int_{f_t(y)}\al + \int_{f_t(x)}\al\\
&=& \int_0^1\int_0^1 d\al(df_t(\dot{\ga}(s)),\dot{f_t}(\ga(s)))dtds 
 - \int_{f_t(y)}\al + \int_{f_t(x)}\al\\
&=& \int_0^1 \int_0^1 -dF_t(df_t(\dot{\ga}(s)))dtds 
 - \int_{f_t(y)}\al + \int_{f_t(x)}\al\\
&=& - \int_0^1 (F_t(f_t(y)) - F_t(f_t(x)))dt 
 - \int_{f_t(y)}\al + \int_{f_t(x)}\al\\
&=& \int_{f_t(x)}\al + \int_0^1F_t(f_t(x))dt
- \left(\int_{f_t(y)} \al + \int_0^1F_t(f_t(y))dt\right)
\end{eqnarray*}
\end{proof}

\begin{remark}\label{R:action}
If $x\in M$ is a fixed point of $f$ and the orbit $f_t(x)$ is
contractible then $\C K_{\al}(f)$ is equal to the value of the
action functional $\C A(F,x)$ up to the above constant 
(see McDuff-Salamon \cite[Section 9]{MR2045629} for the background about
the action functional).
\end{remark}

\subsection{Normalization}\label{SS:normalization}
The restriction of the above cocycle to the subgroup of
compactly supported Hamiltonian diffeomorphisms takes values in the
space of smooth functions, after choosing $\C K_{\al}(f)$
so that it has  compact support. To see that this is possible
we first observe that the function $\C K_{\al}(f)$ is locally constant 
outside the support of $f.$ Let $x,y\in M$ be two points outside
the support of the Hamiltonian function $F_t$ generating $f.$
Notice that,
by definition, a compactly supported Hamiltonian diffeomorphism
is generated by a compactly supported Hamiltonian function.
According to the formula in Proposition \ref{P:action}, we get
that $\C K_{\al}(f)(x)=\C K_{\al}(f)(y)$ and hence $\C K_{\al}(f)$
is actually constant outside the support of $f$ so we can normalize
it to be compactly supported.

Since the integration with respect to the symplectic volume is
a morphism of $\Ham_c(M,d\alpha)$-modules
$C^\infty_c(M)\to\B R,$
the integral integral $\int_M\C K_{\al}(f)\,(d\al)^n$ for normalized
$\C K_{\al}(f)$ is a cocycle in a trivial module thus a homomorphism
$$
\operatorname{Cal}\colon\Ham_c(M,d\al)\to \B R
$$
called the Calabi homomorphism. Its value can be calculated
using the following formula
$$\int_M \C K_{\al}(f) \om^n = (n+1)\int_0^1\int_M F_t\om^n\, dt$$ 
(see Arnold-Khesin \cite[Theorem 8.7 on page 247]{MR1612569}).
Here $F_t$ is any Hamiltonian function generating $f$.

On the other hand, if we only assume that a symplectomorphism has a
compact support, but not that it is generated by a (compactly supported)
Hamitonian flow then some Dehn twists $\tau$ of $T^*S^1$ supported near
the zero section induce $\C K_{pdq}(\tau)$ which cannot be
normalized to have a compact support. More precisely,
let $\tau\colon T^*S^1\to T^*S^1$ be given by
$\tau(p,q)=(p,q+t(p)),$ where $0\leq t(p)\leq 2\pi$,
$t(p)=0$ for $p\leq-1$ and
$t(p)=2\pi$ for $p\geq1$. 
Then $d\C K_{pdq}(\tau)=t'(p)pdp$ and 
\begin{eqnarray*}
\C K_{pdq}(\tau)(1,q)-\C K_{pdq}(\tau)(-1,q)&=&\int_{-1}^1t'(p)pdp\\
&=&t(p)p\, {\text{\Big |}} _{-1}^1-\int_{-1}^1t(p)dp\\
&=&2\pi - \int_{-1}^1t(p)dp.
\end{eqnarray*}
Thus $\C K_{pdq}(\tau)$ has compact support if and only if
$\int_{-1}^1 t(p)dp =2\pi.$

\subsection{The case of a closed manifold}\label{SS:closed}
Let $p\colon\widetilde M \to M$ be the universal covering.
Let  $G$ be a path connected group of homeomorphisms of $M$ and let
$$
\text{ev}_x\colon G \to M
$$
be the evaluation map defined by
$$\text{ev}_x(f):= f(x),$$
where $x\in M$ is a chosen point.

\begin{proposition}[{cf. \cite[Remark 1.5.B]{MR97g:58021}}]\label{P:eval}
%\begin{proposition}\label{P:eval}
Let $G$ be a path connected group of homeomorphisms of $M$.
Suppose that the evaluation map induces the trivial homomorphism on the 
fundamental group. Then there exists a homomorphism
$$\Psi\colon G \to \operatorname{Homeo}(\widetilde M)$$ 
such that
$\Psi(f)\circ \pi = \Psi(f)$, where $\pi\in\pi_1(M)$
is a deck transformation. Moreover it satisfies
$p\circ \Psi(f)= f\circ p$.
\end{proposition}

\begin{proof}
Let $f_t\in G$ be an isotopy from the identity to $f$.
Define
$$
\Psi(f)(x):= \tilde f_1(x),
$$
where $\tilde f_t$ is the lift of the path $f_t(p(x))$
starting at $x\in \widetilde M.$ 

We have to show that $\Psi(f)$ is well defined.
Let $f'_t\in G\Mo$ be another isotopy from the
identity to $f$. Then we get that
$\tilde f_1(x) = \pi(\tilde f'_1(x)),$ where
$\pi:=[f_t(x)*\overline{f'_t(x)}]\in \pi_1(M)$ is
the deck transformation induced by the loop
$(f_t * \overline{f'_t})(x)=ev_x(f_t * \bar{f'_t})$. 
Since the evaluation map induces the trivial homomorphism
on the fundamental groups, the above element $\pi$
is trivial. Hence $\Psi(f)(x)$ is well defined and
clearly it is a lift of $f$.
\end{proof}

\begin{example}\label{E:eval}
\hfill
\begin{enumerate}
\item
It is easy to prove that the image of the homomorphism
$$(\text{ev}_x)_*\colon \pi_1(G)\to \pi_1(M)$$
is contained in the centre of the fundamental group of $M$.
Hence, the above proposition applies to any space whose
fundamental group has trivial centre.
\item
If $G\subseteq\Symp_0\Mo$ and the evaluation map 
$\text{ev}_x\colon G\to M$ induces the
trivial map on the fundamental
then the homomorphism $\Psi$ takes
values in the group $\Symp_0(\widetilde M,p^*\om).$
\end{enumerate}
\end{example}

It is well known (and difficult to prove, see Corollary 9.1.2 in 
McDuff-Salamon \cite{MR2045629}) that the
evaluation map on the group of Hamiltonian diffeomorphisms
of a closed manifold induces the trivial map on the
fundamental group. Hence we get the lifing homomorphism
$\Psi\colon \Ham\Mo\to \Ham(\widetilde M,p^*\om).$
If $\Mo$ is closed symplectically aspherical then
taking the composition $\C K_{\al}\circ \Psi$ we get a cocyle
$$\Psi^*(\C K_{\al})\colon  \Ham\Mo \to C^{\infty}(\widetilde M)/\B R.$$

It follows from the next proposition that for symplectically 
hyperbolic manifolds
the above cocycle is defined on the group of all symplectomorphisms
isotopic with the identity.

\begin{proposition}\label{P:eval_hip}
Let $\Mo$ be a closed symplectically hyperbolic manifold.
Then the evaluation map $\text{ev}_x\colon\Symp_0\Mo\to M$
induces the trivial homomorphism on the fundamental group.
Consequently, there exists a lifting homomorphism
$\Psi\colon\Symp_0(M,\om)\to \Symp_0(\widetilde M,p^*\om).$
\end{proposition}

\begin{proof}
If $\Mo$ is symplectically hyperbolic then the symplectic form
vanishes on tori. That is $s^*[\om]=0$ for any smooth map
$s\colon T^2\to M$. The proof of this fact is not difficult and
can be found in Gromov \cite[Example 0.2.A']{MR1085144} or
K\k edra \cite[Proposition 1.9]{hyperbolic}. This implies that
the flux homomorphism 
$\text{Flux}\colon \pi_1(\Symp_0\Mo)\to H^1(M;\B R)$ is trivial.
Consequently the inclusion $\Ham\Mo\to \Symp_0\Mo$ is a homotopy
equivalence and the statement follows form the analogous
statement for $\Ham\Mo.$
\end{proof}

\begin{corollary}\label{C:hipeval}
If $\Mo$ is a closed symplectically hyperbolic manifold whose
fundamental group has torsion-free centre then there
is a well defined cocycle
$$
\Psi^*(\C K_{\al})\colon\Symp_0\Mo \to C^{\infty}(\widetilde M)/\B R.
$$
\qed
\end{corollary}

%%%%%%%%%%%%%%%%%%%%%%%%%%%%%%%%%%%%%%%%%%%%%%%%%%%%%%%%%%%%%%%
\section{The Polterovich homomorphism}\label{S:polterovich}
%%%%%%%%%%%%%%%%%%%%%%%%%%%%%%%%%%%%%%%%%%%%%%%%%%%%%%%%%%%%%%%

Let $\Ham(M,\om,x,y)$ denote the group of Hamiltonian diffeomorphisms
fixing two points $x,y\in M$. 
Let $\ga\colon [0,1]\to M$ be a curve from $x$ to $y$. That is $\ga(0)=x$ and 
$\ga(1)=y.$ Let $f\in \Ham(M,\om,x,y)$. Let $D\subset M$
be a disc whose boundary is the loop given by the concatenation
of $\ga$ and $f(\ga)$. In other words, let $D\colon D^2\to M$
be a 2-chain with $\partial D=\ga - f(\ga).$ Define the map
$$
\C P_{x,y}\colon \Ham(M,\om,x,y)\to \B R
$$
$$
\C P_{x,y}(f) := \int_D\om.
$$

\begin{proposition}\label{P:polterovich}
Let $\Mo$ be a symplectically aspherical manifold. 
\begin{enumerate}
\item
If $\Mo$ is exact then
the map $\C P_{x,y}$ is a well defined homomorphism.
\item
If $M$ is closed then 
the map $\C P_{x,y}$ is a well defined homomorphism
on the group $\Ham(M,\om,x,y)_{\text{cont}}$ consisting
of those Hamiltonian diffeomorphisms for which both $x$ and
$y$ are contractible fixed points.
\end{enumerate}
In both cases we shall call it the {\bf Polterovich homomorphism}
(cf. Section 2.1 in Polterovich \cite{MR2003i:53126}).
\end{proposition}

A fixed point $x\in M$ of $f$ is called {\bf contractible}
if there exists an isotopy $f_t$ from the identity to $f$ such that
the loop $f_t(x)$ is contractible in $M.$ 

\subsection{A digression on contractible fixed points}
There is an obvious inclusion 
$\Ham(M,\om,x,y)_0\subset \Ham(M,\om,x,y)_{\text{cont}}.$ 
If $\dim M\geq 4$ the above inclusion is an equality.
To see this, consider the evaluation fibration
$$
\Ham(M,\om,x,y)\to \Ham\Mo \stackrel{\text{ev}_{x,y}}\longrightarrow
 M\times M - \Delta,
$$
where $\text{ev}_{x,y}(f) = (f(x),f(y)).$
If both $f_t(x)$ and $f_t(y)$ are contractible in $M$ then
the loop $(f_t(x),f_t(y))$ is contractible in $M\times M - \Delta$ because
the diagonal has codimension at least four. Lifting the contraction
from $(f_t(x),f_t(y))$ to $(x,y)$ we get an isotopy 
from the identity to $f$ fixing $x$ and $y$ as claimed.

If $\dim M=2$ then it is easy to construct a Hamiltonian diffeomorphism $f$  
of a surface with two contractible fixed points such that any isotopy
from the identity to $f$ does not preserve them.

\subsection{Back to the proof of Proposition \ref{P:polterovich}}

We start with the following lemma which is of independent interest.
Let $x,y\in M$ and let $(x-y)\colon C^{\infty}(M)/\B R\to \B R$
denote the evaluation homomorphism. That is $(x-y)(F):= F(x) - F(y)$.

\begin{lemma}\label{L:KP}
Let $(M,d\al)$ be an exact symplectic manifold.
The following diagram is commutative.
$$
\xymatrix
{
\Ham(M,d\al) \ar[r]^{\C K_{\al}} & C^{\infty}(M)/\B R \ar[d]^{x-y}\\
\Ham(M,d\al,x,y) \ar[r]^{\phantom{padupad}\C P_{x,y}} \ar[u] & \B R\\
}
$$
\end{lemma}

\begin{proof}
This is an application of the Stokes Lemma. Let $D$ be a singular
disc as in the previous section.
\begin{eqnarray*}
\C P_{x,y}(f)&=& \int_D d\alpha 
= \int_{\partial D}\alpha\\
&=& \int_{\ga} \alpha - \int_{f(\ga)}\alpha
= \int_{\ga} \alpha - f^*\alpha\\
&=& -\int_{\ga}d\C K_{\al}(f)
= -\int_{\partial \ga}\C K_{\al}(f)\\
&=& \C K_{\al}(f)(x) - \C K_{\al}(f)(y).
\end{eqnarray*}
\end{proof}

\begin{proof}[Proof of Proposition \ref{P:polterovich}]

Since $(x-y)\colon C^{\infty}(M)/\B R \to \B R$ is a morphism
of  $\Ham(M,\om,x,y)$-modules, the composition
$$(x-y)\circ K_{\al}\colon \Ham(M,d\al,x,y)\to \B R$$ 
is a cocycle in the trivial module, hence a homomorphism.
This proves the proposition in the case of an exact manifold.

In the closed case, let $\tilde x,\tilde y\in \widetilde M$ be
points in the preimage $p^{-1}(x)$ and $p^{-1}(y)$ respectively.
As in Proposition \ref{P:eval} and due to the contractibility
of the fixed points, we have a homomorphism
$$
\Psi\colon \Ham(M,\om,x,y)_{\text{cont}}\to 
\Symp(\widetilde M,d\al,\tilde x,\tilde y).
$$
Clearly, the composition $(\tilde x-\tilde y)\circ \C K_{\al}\circ \Psi$
is a homomorphism equal to $\C P_{x,y}$.
\end{proof}

%%%%%%%%%%%%%%%%%%%%%%%%%%%%%%%%%%%%%%%%%%%%%%%%%%%%%%%%%%%%%%%%%%%%%%%%%%%%%
\section{The non-triviality of $\C K$}\label{S:nontrivial}
%%%%%%%%%%%%%%%%%%%%%%%%%%%%%%%%%%%%%%%%%%%%%%%%%%%%%%%%%%%%%%%%%%%%%%%%%%%%%

\begin{theorem}[Frauenfelder-Polterovich-Schlenk-Schwarz]\label{T:PS}
Let $\Mo$ be a symplectic manifold which is either 
\begin{enumerate}
\item
closed symplectically aspherical or
\item
exact and convex.
\end{enumerate}
Let $f\in \Ham_c(M,\om)$ be a compactly supported
Hamiltonian diffeomorphism. 
There exist two contractible points
$x,y\in M$ such that $\C P_{x,y}(f)\neq ~0.$
\end{theorem}

\begin{proof}
First observe that, according to Remark \ref{R:action}, 
we have
$$
\C P_{x,y}(f)=\C K_{\al}(f)(x) - \C K_{\al}(f)(y)
=\C A(F,x) - \C A(F,y),
$$
for any two contractible fixed points $x,y\in M.$ 
This observation has been made by Polterovich in
\cite{MR2003i:53126}.

The second step is to prove that for any Hamiltonian
diffeomorphism $f\in \Ham\Mo$ there exist two contractible
fixed points $x,y\in M$ such that $\C A(F,x) < \C A(F,y).$
In the case of a closed symplectically aspherical manifold
this was done by Schwarz in \cite{MR1755825}. We also refer to
McDuff-Salamon (Theorem 9.1.6 on page 283 in \cite{MR2045629})
for the proof. 

The case of an exact and convex symplectic manifold was done
by Frauenfelder and Schlenk \cite{MR2342472} based on the ideas
of Schwarz. More precisely they define (Section 7 of \cite{MR2342472}) 
a map 
$$c\colon \Ham_c\Mo\to \B R$$ 
and they prove (Theorem 7.3 in \cite{MR2342472}) that for any
compactly supported Hamiltonian diffeomorphism $f$
the number (called the Schwarz norm of $f$) 
$$c(f)+c(f^{-1})$$ 
is positive if and only if $f\neq \text{ID}$. Moreover, it follows from 
Proposition 5.1 and Corollary 6.3 in \cite{MR2342472}
that there exist two fixed points $x,y\in M$ contractible
with respect to the isotopy $f_t$ generated by $F$ such that
$c(f)=\C A(F,x)$ and $c(f^{-1}) = -\C A(F,y)$.
Finally we get that
$$
0<c(f)+c(f^{-1})=\C A(F,x) - \C A(F,y) = \C P_{x,y}(f),
$$
as claimed.
\end{proof}

\begin{corollary}\label{C:nontrivial_exact}
Let $\Mo$ be an exact and  symplectic manifold convex at infinity.
Then the cohomology class 
$[\C K_{\al}]\in H^1(\Ham_c(M,\om);C_c^{\infty}(M))$
is non-zero.
\end{corollary}

\begin{proof}
Let $f\in \Ham_c(M,\om,x,y)$ be such that
$\C P_{x,y}(f)\neq 0$. 
Let 
$$i\colon \Ham_c(M,\om,x,y)\to \Ham_c(M,\om)$$
be the inclusion. Consider the composition

$$
\xymatrix
{
H^1(\Ham_c(M,\om);C_c^{\infty}(M)) 
\ar[d]^{i^*} & \\
H^1(\Ham_c(M,\om,x,y);C_c^{\infty}(M)) 
\ar[r]^{\phantom{dupa}(x-y)_*} & 
H^1(\Ham_c(M,\om,x,y);\B R),
}
$$
where the second map is induced by the morphism 
$x-y\colon C_c^{\infty}(M)\to \B R$
of $\Ham_c(M,\om,x,y)$-modules. 

Notice that $\B R$ is the trivial module hence
$H^1(\Ham_c(M,\om,x,y);\B R)=
\Hom(\Ham_c(M,\om,x,y),\B R).$
Clearly the class $(x-y)_*i^*[\C K]= [\C P_{x,y}]$
is non-trivial since $\C P_{x,y}(f)\neq 0$.
Thus the class $[\C K]\neq 0$ 
%$H^1(\Ham_c(\B R^{2n},\om);C_c^{\infty}(\B R^{2n}))$
as claimed.
\end{proof}

An analogous proof works in the closed case, so we get the following.

\begin{corollary}\label{C:nontrivial}
Let $\Mo$ be a closed symplectically aspherical manifold.
Then the cohomology class 
$\Psi^*[\C K_{\al}]\in H^1(\Ham(M,\om);C^{\infty}(\widetilde M)/\B R)$
is non-zero.\qed
\end{corollary}

%%%%%%%%%%%%%%%%%%%%%%%%%%%%%%%%%%%%%%%%%%%%%%%%%%%%%%%%%%%%%%%%%%
\section{Distortion of cyclic subgroups}\label{S:undistorted}
%%%%%%%%%%%%%%%%%%%%%%%%%%%%%%%%%%%%%%%%%%%%%%%%%%%%%%%%%%%%%%%%%%

\subsection{Abstract nonsense}\label{SS:abstract}
Let $\Gamma $ be a finitely generated group equipped with
a word norm $|\cdot|$ and let $(V,\|\cdot\|)$ be a $\Gamma $-module
with a norm preserved by the action. Let $\left<f\right>\subset \Gamma $
be a cyclic group generated by $\text{Id}\neq f\in \Gamma $.
Consider the following commutative diagram 
$$
\xymatrix
{
\Gamma  \ar[r]^{\C K} & V \ar[d]^{\mu}\ar[r]^{\|\cdot\|} & \B R\\
\left<f\right> \ar[u] \ar[r]^{\C P} & \B R
}
$$
where the maps $\C K$ and $\C P$  are cocycles, 
$\mu\colon V \to \B R$ is a map of $\left<f\right>$-modules and
$\B R$ is a trivial $\left<f\right>$-module.

\begin{lemma}\label{L:nonsense}
Suppose that the following conditions hold:
\begin{enumerate}
\item
$|\C P(f^k)|\leq \|\C K(f^k)\|$
for all $k\in \B Z$, 
\item
$\C P$ is a nontrivial homomorphism,
\end{enumerate}
Then $\left<f\right>$ is undistorted in $\Gamma $.
\end{lemma}

\begin{proof}
First observe that any cocycle $\C K\colon \Gamma  \to V$ is
Lipschitz. To see this, let $g_1,\ldots,g_k\in \Gamma $ be
the generators inducing the norm $|\cdot|$. 
Let $m:=\max\{\|\C K(g_i)\| \,|\, i=1,\ldots,k\}.$
If $f=g_{i_1}\cdot \ldots \cdot g_{i_{|f|}}$
then we get that
\begin{eqnarray*}
\|\C K(f)\| &=& \|\C K(g_{i_1}\cdot \ldots \cdot g_{i_{|f|}})\|\\
&=&\left \|\sum_{j=1}^{|f|}\C K(g_{i_j})\cdot g_{i_{j+1}}\cdot \ldots \cdot g_{i_{|f|}}\right \|\\
&\leq& \sum_{j=1}^{|f|}\|\C K(g_{i_j})\| \leq m\cdot |f|.
\end{eqnarray*}

Let $\text{Id}\neq f\in \Gamma $. The hypothesis together with the above
observation immediately imply the following inequalities:
$$
0<n|\C P(f)|=|\C P(f^n)| \leq \|\C K(f^n)\| \leq m\cdot |f^n|.
$$
It follows that
$$
|f^n| \geq \frac{|\C P(f)|}{m} \cdot n>0,
$$
which finishes the proof.
\end{proof}

\subsection{Proof of Theorem \ref{T:undistorted}}\label{SS:undistorted}
Recall that $(M,d\al)$ is an exact and convex symplectic manifold and
$\Gamma \to \Ham_c(M,d\al)$ is an effective action of a finitely generated group.
Let $\text{Id}\neq f\in \Gamma $. We shall show that the cyclic subgroup
$\left<f\right>\subset \Gamma $ is undistorted.

Let the {\bf oscillation} $\text{osc}\colon C^{\infty}_c(M)\to \B R$ be defined by
$\text{osc}(F):= \max_{x\in M} F - \min_{x\in M} F$. It is a norm
on the space of compactly supported functions.
We have the following commutative diagram
$$
\xymatrix
{
\Gamma  \ar[r] & \Ham_c(M,d\al) \ar[r]^{\phantom{dupa}\C K_{\al}} &
C_c^{\infty}(M)\ar[d]^{x-y} \ar[r]^{\phantom{dup}\text{osc}} & \B R\\
\left<f\right>\ar[u]\ar[r] & \Ham(M,d\al,x,y) \ar[u]\ar[r]^{\phantom{dupadupa}\C P_{x,y}} & \B R
}
$$
Observe that the inequality $\C P_{x,y}(f) \leq \text{osc}\C K_{\al}(f)$ 
is trivially satisfied since 
$\C P_{x,y}(f) = \C K_{\al}(f)(x) - \C K_{\al}(f)(y)$, according to 
Lemma \ref{L:KP}. The points $x,y \in M$ are chosen so that
the Polterovich homomorphism is nontrivial $\C P_{x,y}(f)>0.$
Hence we are precisely in the situation when the hypothesis of 
Lemma \ref{L:nonsense} is satisfied. Thus applying this lemma 
finishes the proof.
\qed

%%%%%%%%%%%%%%%%%%%%%%%%%%%%%%%%%%%%%%%%%%%%%%%%%%%%%%%%%%%%%%%%%%%%%%%%%
\subsection{An example of a compactly supported action with 
distorted cyclic subgroups}
\label{SS:heisenberg}
%%%%%%%%%%%%%%%%%%%%%%%%%%%%%%%%%%%%%%%%%%%%%%%%%%%%%%%%%%%%%%%%%%%%%%%%%%
This section is essentially due to \L ukasz Gra\-bow\-ski.
If $\Gamma$ acts on a manifold $M$ with a global fixed point $x$
and the isotropy representation can be deformed to the trivial
one. This means that the isotropy representation lies in the
connected component of the trivial one inside the representation variety
$\{\Gamma\to\operatorname{GL}(T_xM)\}$.
%there exist a continous path joining the isotropy
%representation $\Gamma\to \text{GL}(T_xM)$ with the
%trivial ($m$-dimesional) one.
Then one can replace $x$ by a disc and extend the action of $\Gamma$
to an action which is trivial in some (smaller) disc.
In other words one can produce a compactly supported action of
$\Gamma$ on $M-\{x\}$ which has essentially the same orbit structure.

Consider the projective action of the Heisenberg group $H:=N(3,\B Z)$ of
uppertriangular matrices with integer entries on two dimensional sphere $\B S^2$.
The isotropy representation on the tangent space to the north pole $(1,0,0)$ is
given by
$$
\left [
\begin{array}{ccc}
1 & a & c\\
0 & 1 & b\\
0 & 0 & 1\\
\end{array}
\right ] \mapsto
\left [
\begin{array}{cc}
1 & b \\
0 & 1 \\
\end{array}
\right ]
$$
and it clearly can be deformed to the trivial one.
Thus we get that the
Heisenberg group acts smoothly with compact supports on an open ball
$\B S^2 -\{x\}$. This can be generalized in an obvious way
to higher dimensions. The centre of the Heisenberg group
is a distorted cyclic subgroup as it can be directly calculated
(see also Gromov \cite[page 52]{MR1253544}).
We have then proved the following result.

\begin{proposition}\label{P:heisenberg}
There exists a smooth, compactly supported action of the Heisenberg
group $H$ on $\B R^n$ for any $n\geq 2$.\qed
\end{proposition}

%%%%%%%%%%%%%%%%%%%%%%%%%%%%%%%%%%%%%%%%%%%%%%%%%%%%%%%%%%%%%%%%%%%%%%%
\section{A proof of the Polterovich theorem}\label{S:osc}
%%%%%%%%%%%%%%%%%%%%%%%%%%%%%%%%%%%%%%%%%%%%%%%%%%%%%%%%%%%%%%%%%%%%%%%

\begin{proposition}\label{P:osc}
Let $\Mo$ be a closed symplectically hyperbolic manifold.
Then the composition
$$
\Ham\Mo \stackrel{\Psi}\to \Ham(\widetilde M,d\al)
\stackrel{\C K_{\al}}\to C^{\infty}(\widetilde M)/\B R
$$
takes values in $C_b^{\infty}(\widetilde M)/\B R$,
where $C_b^{\infty}(\widetilde M)$ denotes the space of
bounded functions.
\end{proposition}

\begin{proof}
Let $f\in\Ham\Mo$ and
let $x,y\in \widetilde M$ be any two points. Let 
$f_t$ be an isotopy form the identity to $f$ generated by
a Hamiltonian function $F_t\colon M\to \B R$ and let
$\tilde f_t$ and $\tilde F_t$ be the corresponding lifts
to $\widetilde M$.
According to
Proposition \ref{P:action} we get the first equality in the
following calculation.

\begin{eqnarray*}
&&|\,\C K_{\al}(\Psi(f))(x) - \C K_{\al}(\Psi(f))(y)\,|\\
&=&\left |\,\int_{\tilde f_t(x)}\al + \int_0^1\tilde F_t(\tilde f_t(x))dt 
- \int_{\tilde f_t(y)}\al + \int_0^1\tilde F_t(\tilde f_t(y))dt\, \right |\\
&\leq & 2\cdot C \cdot \max_{x}\text{Length}(\tilde f_t(x)) +
2\cdot \max_{x,t} \tilde F_t(x) 
\end{eqnarray*}
The last quantity is finite because $\tilde F_t$ and $\tilde f_t$
are lifts of $F_t$ and $f_t$ respectively and the latter
are defined on a compact manifold. Also, the length is calculated
with respect to the metric induced from $M$.
We used the easy fact that
$$
\int_{f_t(x)}\al \leq C\cdot \text{Length}(f_t(x)),
$$
where $C>0$ is the constant bounding the norm of $\al$.
\end{proof}

Let $\Gamma \to \Ham\Mo$ be an effective action of a finitely generated
group. We shall show that every nontrivial cyclic subgroup
$\left<f\right>$ of $\Gamma $ is undistorted.

\begin{proof}[Proof of Theorem \ref{T:polterovich}]
Let $g_1,g_2,\ldots,g_k\in ~\Gamma $ be generators and let
$m:=\max_i \C K_{\al}(\Psi(g_i))$ be the Lipschitz constant (cf. Lemma \ref{L:nonsense}). 
Due to the above proposition,
$m$ is finite and we have the following estimate
$$
|f^n| \geq \frac{\text{osc}(\C K_{\al}(\Psi(f^n)))}{m}\geq 
\frac{\C P_{x,y}(f^n)}{m} = n \,\frac{\C P_{x,y}(f)}{m}>0.
$$
The existence of fixed points $x,y\in M$ such that
$\C P_{x,y}>0$ is ensured by Theorem \ref{T:PS}.
\end{proof}

%%%%%%%%%%%%%%%%%%%%%%%%%%%%%%%%%%%%%%%%%%%%%%%%%%%%%%%%%%%%%%%%
\section{The relation with the flux homomorphism}\label{S:flux}
%%%%%%%%%%%%%%%%%%%%%%%%%%%%%%%%%%%%%%%%%%%%%%%%%%%%%%%%%%%%%%%%
Let $p\colon\widetilde M \to M$ be the universal cover of a closed manifold. 
Let $C_b^{\infty}(\widetilde M)$ denote the space of bounded 
functions and let
$$\Phi\colon H^1(M;\B R)\to C^{\infty}(\widetilde M)/C_b^{\infty}(\widetilde M)$$
be a homomorphism defined by $d\Phi[a] = p^*a.$ It is easy to check that
this homomorphism is injective.

\begin{proposition}\label{P:flux}
Let $\Mo$ be a closed symplectically hyperbolic manifold.
The following diagram is commutative.
$$
\xymatrix
{
%\Ham\Mo \ar[d]\ar[r] & 
\Symp_0\Mo \ar[d]^{\Psi}\ar[r]^{\text{Flux}} &H^1(M;\B R)\ar[d]^{\Phi}\\
%\Symp_b(\widetilde M,d\al)\ar[r] & 
\Symp(\widetilde M,d\al)\ar[r]^{\C F_{\al}}
&C^{\infty}(\widetilde M)/C_b^{\infty}(\widetilde M)
}
$$
Here $\C F_{\al}$ is the composition of $\C K_{\al}$ with the natural projection.
\end{proposition}

\begin{proof}
We shall show that $\Phi(\Flux(f))$ equals $\C K_{\al}(\Psi(f))$ up
to a bounded summand. 
Notice that the homomorphism $\Psi\colon\Symp_0\Mo\to \Symp_0(\widetilde M,p^*\om)$
exists due to Proposition \ref{P:eval_hip}.
If $f_t$ is an isotopy from the identity to $f$ then
$$\int_0^1(\iota_{\dot{f_t}}\om)\,dt\in \Om^1(M)$$ 
is a differential form representing the flux of $f$. 
In the following calculation let $\C X_t:= \frac{d}{dt}\Psi(f_t).$
\begin{eqnarray*}
d\Phi(\Flux(f))&=&
p^*\left(\int_0^1(\iota_{\dot{f_t}}\om)\,dt\right) \\
&=& \int_0^1(\iota_{\C X_t}d\al)\,dt\\
&=& \int_0^1(\iota_{\C X_t}\Psi(f_t)^*d\al)\,dt\\
&=& \int_0^1\left(\C L_{\C X_t}\Psi(f_t)^*\al - d\iota_{\C X_t}\Psi(f_t)^*\al\right)\,dt\\
&=& \Psi(f)^*\al - \al  - d\int_0^1\Psi(f_t)^*\al(\C X_t)\,dt\\
&=& d\left(\C K_{\al}(\Psi(f)) - \int_0^1\Psi(f_t)^*\al(\C X_t)\,dt\right)\\
\end{eqnarray*}
Further we have that
$$
\left|\int_0^1\Psi(f_t)^*\al(\C X_t)\,dt\right |
\leq C\cdot\max_{(x,t)\in M\times [0,1]}\|D_x\Psi(f_t)\|\cdot \text{Length}(f_t(x)),$$
where $C$ is a constant bounding the norm of the primitive $\al$.
Clearly, the right hand side of the above inequality is finite.
This finishes the proof.
\end{proof}

As a byproduct of the above result we get the following 
characterization of Hamiltonian diffeomorphisms of a closed
symplectically hyperbolic manifold.

\begin{corollary}\label{C:flux}
Let $\Mo$ be a closed symplectically hyperbolic manifold.
A symplectic diffeomorphism $f\in \Symp_0\Mo$ is Hamiltonian
if and only if  $\C K_{\al}(\Psi(f))$ is bounded.\qed
\end{corollary}

%\nocite{*}
\bibliography{../../bib/bibliography}
\bibliographystyle{smfplain}
\end{document}